\documentclass[12pt, reqno]{amsart}

\usepackage{amsfonts,amsmath,amsthm,amssymb}
\usepackage{latexsym}

\begin{document}

{\theoremstyle{plain}
  \newtheorem{theorem}{Theorem}[section]
  \newtheorem{corollary}[theorem]{Corollary}
  \newtheorem{proposition}[theorem]{Proposition}
  \newtheorem{lemma}[theorem]{Lemma}
  \newtheorem{question}[theorem]{Question}
  \newtheorem{conjecture}[theorem]{Conjecture}
}

{\theoremstyle{definition}
  \newtheorem{definition}[theorem]{Definition}
  \newtheorem{remark}[theorem]{Remark}
  \newtheorem{example}[theorem]{Example}
}

\numberwithin{equation}{section}

\bibliography{file 1, file 2, ....}
\bibliographystyle{style name}\bibliography{file 1, file 2, ....}

\begin{thebibliography}{10}

\bibitem{B} A. Brown, A structure theorem for a class of grade three perfect ideals, Ph. D. dissertation, {\it Brandies Univ. Waltham, Mass} May (1984).

\bibitem{B-H} W.Bruns, J. Herzog, Cohen-Macaulay rings, {\it Cambridge Stud. Adv. Math., } vol. \textbf{39}, Cambridge University Press, Cambridge (1993).


\bibitem{B-E} D. Buchbaum, D. Eisenbud, Algebra structures for finite free resolutions, and some structure theorems for ideals of codimension 3. {\it Amer. J. Math., } vol. \textbf{99}, No. 3447-485 (1977).

\bibitem{D} S. J. Diesel, Irreduciblility and dimension theorems for families of height 3 Gorenstein algebras. {\it Pacific J. Algebra, } vol. \textbf{172}, No. {\bf 2} 365-397 (1996).

\bibitem{E} D. Eisenbud, Commutative Algebra with a view toward Algebraic Geometry. {\it Grad. texts in Math., } vol. \textbf{150}, Springer (1995).

\bibitem{E} B. Engheta, Bound on the multiplicity of almost complete intersections, preprint.

\bibitem{G-S-S} L. Gold, H. Schenck, H, Srinivasan, Betti numbers and degree bounds for some linked zero-schemes {\it J. Pure Appl. Algebra,  } vol. \textbf{210}, No. {\bf 2} 481-491 (2007).

\bibitem{H-S} J. Herzog, H. Srinivasan, Bounds for multiplicities, {\it Trans. Amer. Math. Soc.,} \textbf{350} 2879-2902 (1998).

\bibitem{H-M} C. Huneke, M. Miller, A note on the multiplcity of Cohen-Macaulay algebras with pure resolutions, {\it Can. J. Math.,} \textbf{37} 1149-1162 (1985).

\bibitem{M-N-R} J. Migliore, U. Nagel, T. R\"omer, The multiplicity conjecture in low codimension, {\it Math. Res. Lett., } vol. \textbf{12} 731-747 (2005).

\bibitem{M} R. M. Mir\'{o}-Roig, A note on the multiplicity of determinantal ideals, {\it J. Algebra,  } vol. \textbf{299}, No. {\bf 2} 714-724 (2006).

\bibitem{P-S} C. Peskine, L. Szpiro, Liaison des vari\'{e}t\'{e}s alg{\'{e}}briques. I, {\it Invent. Math.,} \textbf{26} 271-302 (1974).

\end{thebibliography}

\bibliographystyle{style name}



\title{MULTIPLICITY OF CODIMENSION THREE ALMOST COMPLETE INTERSECTIONS}

\author{
Sumi Seo and Hema Srinivasan}

\address{address: Hema Srinivasan, Depart. of Mathematics, University of Missouri-Columbia}
\email{hema@math.missouri.edu}

\address{address: Sumi Seo, Depart. of Mathematics, University of Missouri-Columbia}
\email{sumi@math.missouri.edu}

\thanks{}

\subjclass{}

\begin{abstract}
 We establish the upper bound in the multiplicity conjecture of Herzog and Srinivasan for the codimension three almost complete intersections.  The proof is essentially by direct computation and uses the structure theorem of codimension three almost complete intersection of Buchsbaum and Eisenbud.  We also give some partial results in the case when I is the almost complete intersection ideal linked to a complete intersection in one step. 
\end{abstract}

\maketitle



\newcommand\sfrac[2]{{#1/#2}}

\newcommand\cont{\operatorname{cont}}
\newcommand\diff{\operatorname{diff}}


\section{Introduction}%

Let $R$ be a standard graded polynomial ring $k[x_1, x_2, \cdots x_n]$ where k is a field and let $I\subset R$ be a homgeneous ideal of codimension $h$.

Consider the graded minimal free $R$-resolution of $R/I$:

\[
0\,\longrightarrow\, \oplus_{j\in\mathbb Z}R(-j)^{\beta_{p,j}}
\,\longrightarrow\,\ldots \,\longrightarrow\, \oplus_{j\in\mathbb
Z}R(-j)^{\beta_{1,j}}
\,\longrightarrow\,R\,\longrightarrow\,R/I\,\longrightarrow\, 0
\]\\
where $\beta_{i,j} = dim Tor_{i}^{R}(R/I,k)_{j}$ the
graded Betti number of $R/I$ and $p$ is the projective dimension of
$R/I$.

Let $h$ denote the codimension of $R/I$ and let $e(R/I)$ be the
multiplicity of $R/I$. Then $h\le p$ and equality holds if and only
if $R/I$ is Cohen-Macaulay. Let  $m_{i}(I)=min\{j\in \mathbb
Z\mid\beta_{i,j}(R/I)\neq 0\}$  be the minimal 
and $M_{i}(I)=max\{j\in\mathbb Z\mid\beta_{i,j}(R/I)\neq 0\}$ be the
  maximal  shifts at the i-th step.

\textbf{Conjecture 1.1} (Herzog-Huneke-Srinivasan) \ If $R/I$ is
Cohen-Macaulay, then\

\[
\frac{\prod_{i=1}^{h}m_{i}}{h!}\le
e(R/I)\le\frac{\prod_{i=1}^{h}M_{i}}{h!}.
\]

Conjecture 1.1 has been studied in many people and has resulted in
some achievement with some conditions.

\textbf{Conjecture 1.2} (Herzog-Srinivasan) Even if $R/I$ is not
Cohen-Macaulay, the multiplicity  $e(R/I)$ satisfies\

\[
e(R/I)\le\frac{\prod_{i=1}^{h}M_{i}}{h!}.
\]

These conjectured bounds have been established in many cases \footnote { After we submitted this paper we heard the exciting news that the multiplicity conjecture is solved for Cohen Macaulay $R/I$ in  characterisitc zero by Eisenbud and Schreyer. }.
  For instance,
when $R/I$ is Cohen-Macaulay with a pure resolution, that is $m_i = M_i= d_i$ for all $1\le i\le p$, then
the conjecture follows from the formula $e(R/I) = {{\prod _{i=1}^p d_i} \over {p!}}$ in \cite{H-M}.
It is known in codimension two by results of Herzog-Srinivasan, Gold and R\"omer.
I is known in Gorenstein codimension two by the results of Herzog- Srinivasan, Migliore-Nagel-R\"omer.   The very next codimension three case is that of almost complete intersection.   Recently, Bahman Engheta \cite{E}  has also studied these algebras  and has  obtained other  bounds for the multiplicity in terms of the degrees of the generators.

  In this paper, our main theorem proves the upper bound for almost complete intersections of codimension three.  We discuss the lower bound in some cases.  We give some examples and some structure of the resolutions of almost complete intersections in terms of the linked Gorenstein ideal.  The main tool in the proof is the
  theorem of Buchsbaum and Eisenbud which states that all codimension three almost complete intersections arise as $(K:J)$ where $J$ is a codimension three Gorenstein ideal and $K$ is a regular sequence of length three contained in J.   The proof is then
 by direct computations.

\section{The Multiplicity of Almost Complete Intersection with Codim 3}%

Let $ R = k[x_1\cdots x_d] $  and  $I$  be a homogeneous almost complete intersection ideal of codimension three.    By a theorem of Buchsbaum and Eisenbud in \cite{B-E},
there exists a Gorenstein ideal $J$ and a regular sequence $K = (f_1, f_2, f_3)$ contained
in $J$ such that $I = (K:J)$ and $J = (K:I)$.  By Peskine and Szpiro in \cite{P-S}, a resolution of $R/I$ over $R$ can be obtained as the dual of the mapping cone of the resolution of
$R/J$ and the resolution of $R/K$.  By Buchsbaum and Eisenbud structure theorem on Gorenstein ideals, we know that $J$ is minimally generated by an odd number 2n+1 of
elements in $R$ which form the $2n\times 2n$ order pfaffians of a $2n+1 \times 2n+1$ skew symmetric matrix.  In fact, since J is homogeneous this matrix can be taken to be
homogeneous as well \cite{D}.  Since I is homogeneous we can extend the entire argument
to take $f_1, f_2, f_3$ as homogeneous elements in $J$ and the mapping cone of the
dual of the graded resolutions of $R/J$ and $R/K$ induced by the inclusion of $K$ into $J$ will be a graded resolution of $R/I$.  We summarize these as follows:

\begin{theorem}\cite{B-E}  \label{Thm 1}

Suppose that $I=(f_1,f_2,f_3,f_4)$ is a homogeneous
almost complete intersection of codim 3 with $f_1,f_2,f_3$ forming a regular sequence of degree $deg f_i=e_i$.\\
$1$. Then $J=((f_1,f_2,f_3):I)$ is a codimension 3 Gorenstein ideal minimally generated by homogeneous elements
$g_1,\ldots,g_{2m+1}$ of degrees $d_i=deg g_i$ for $i=1,2,\ldots,2m+1$.\\
$2$. Additionally, if $c=\frac{1}{m} \sum_{i=1}^{2m+1}d_i$, $R/J$ has the resolution

\begin{eqnarray} \label{Jres}
0\,\rightarrow\,
R(-c)\,\rightarrow\,\sum_{i=1}^{2m+1}R(-(c-d_{i}))
\,\rightarrow\,\sum_{i=1}^{2m+1}R(-d_{i})
\,\rightarrow\,R\,\rightarrow\,R/J\,\rightarrow\, 0.
\end{eqnarray}

 and $R/K$ has the resolution

\begin{eqnarray}\label{Kres}
0\,\rightarrow\, R(-\sum_{i=1}^3e_i)\,\rightarrow\,\hspace{0.01pt}\sum_{1\le
i<j\le 3}R(-(e_i+e_j)) \,\rightarrow\,\sum_{i=1}^{3}R(-e_{i})
\,\rightarrow\,R\,\rightarrow\,R/K\,\rightarrow\, 0,
\end{eqnarray}

then the resolution of $R/I$ is

\[
0\,\rightarrow\,
\sum_{i=1}^{2m+1}R(-(\sum_{j=1}^3e_j-d_i))\,\rightarrow\,
\sum_{i=1}^{2m+1}R(-(\sum_{j=1}^3e_j-(c-d_i))) \oplus\sum_{1\le
i<j\le 3}R(-(e_i+e_j)) \,
\]

\begin{eqnarray}\label{Ires}
\rightarrow\,R(-(\sum_{j=1}^3e_j-c))\oplus\sum_{j=1}^3R(-e_j)\,\rightarrow\,R\,\rightarrow\,
R/I\rightarrow\, 0.\
\end{eqnarray}
\end{theorem}

The following corollary will be useful in our calculations. We may
assume $e_1\leq e_2\leq e_3$ and $d_1\leq d_2\leq\cdots\leq d_{2m+1}$.

\begin{corollary} \label{Cor1}
Suppose that $I$ is an almost complete intersection of codim 3 with
a regular sequence $(f_1,f_2,f_3)$ among a minimal generating set of
$I$. Suppose that $J=((f_1,f_2,f_3):I)$ is the corresponding
Gorenstein ideal. Then there exists $f_4\in I$ such that $I=(f_1,f_2,f_3,f_4)$ and $deg
f_4=e_1+e_2+e_3-c$ where $mc=\sum_{i=1}^{2m+1}d_i$, $deg f_i = e_i$ and $d_i's$ are the degrees of the minimal generators of $J$. 
\begin{proof}
By the above theorem \ref{Thm 1}, (2.3) is the resolution of $R/I$.
So we can choose the fourth generator for I to be of degree
$e_1+e_2+e_3-c$.
\end{proof}
\end{corollary}

\begin{remark}\cite{P-S} \label{Rem}
The multiplicity $e(R/I)$ can be obtained from the shifts in a graded resolution of $R/I$. If

\[
0\,\rightarrow\, \oplus_{j\in\mathbb Z}R(-j)^{\beta_{p,j}}
\,\rightarrow\,\ldots \,\rightarrow\, \oplus_{j\in\mathbb
Z}R(-j)^{\beta_{1,j}}
\,\rightarrow\,R\,\rightarrow\,R/I\,\rightarrow\, 0
\]

is a graded of $R/I$, then
\begin{equation*}
e(R/I)=(-1)^h \frac{1}{h!} \sum_{j} \sum_{i=0}^{p}(-1)^{i} \beta_{i,j}j^{h}.
\end{equation*}
\end{remark}

\begin{lemma}\label{Lem2}
Suppose that $K=(f_1,f_2,f_3)$ is a regular sequence contained in a codimension three
Gorenstein ideal $J$ generated by $g_1,g_2,\cdots,g_{2m+1}$
.\\ If $I=(K:J)$, then $e(R/I)=e(R/K)-e(R/J)$.
\begin{proof}
Let $degf_i=e_i,\  1 \leq i \leq 3$ and $deg g_j=d_j,\  1 \leq j \leq 2m+1$.
We obtain $e(R/K)=e_1e_2e_3$ since K is a regular sequence and $6e(R/J)=\sum_{i=1}^{2m+1}d_i(c-d_i)(c-2d_i)$ by \cite{H-S}  since $J$ is a Gorenstein ideal. The above remark \ref{Rem} enables us to have
\begin{flalign*}
6e(R/I)&=\sum_{i=1}^{2m+1}(\sum_{j=1}^{3}e_j-d_i)^3-\sum_{i=1}^{2m+1}(\sum_{j=1}^{3}
e_j-(c-d_i))^3\\
&\ \ -((e_1+e_2)^3+(e_1+e_3)^3+(e_2+e_3)^3)+(\sum_{j=1}^{3}
e_j-c)^3+(e_1^3+e_2^3+e_3^3)\\
&=\sum_{i=1}^{2m+1}((\sum_{j=1}^{3} e_j)^3-3(\sum_{j=1}^{3} e_j)^2d_i+3(\sum_{j=1}^{3}
e_j)d_i^2-d_i^3)\\
&\ \  -\sum_{i=1}^{2m+1}((\sum_{j=1}^{3}e_j)^3-3(\sum_{j=1}^{3}
e_j)^2(c-d_i)+3(\sum_{j=1}^{3} e_j)(c-d_i)^2-(c-d_i)^3)\\
&\ \ -((e_1+e_2)^3+(e_1+e_3)^3+(e_2+e_3)^3)+(e_1^3+e_2^3+e_3^3)\\
&\ \ +(e_1+e_2+e_3)^3-3(e_1+e_2+e_3)^2c+3(e_1+e_2+e_3)c^2-c^3\\
&=(e_1+e_2+e_3)^3-((e_1+e_2)^3+(e_1+e_3)^3+(e_2+e_3)^3)+(e_1^3+e_2^3+e_3^3)\\
&\ \ -\sum_{i=1}^{2m+1}d_i(c-d_i)(c-2d_i)\\
&=6e_1e_2e_3-\sum_{i=1}^{2m+1}d_i(c-d_i)(c-2d_i)=6e(R/K)-6e(R/J).
\end{flalign*}
\end{proof}
\end{lemma}
From the structure of the graded resolution (\ref{Jres}) of the Gorenstein ideal
$J$, one can easily see the following numerical criterion on the
degrees of the generators.  This is recorded by
Diesel, in the generalization of Buchsbaum-Eisenbud structure theorem to the graded case \cite{D}.

\begin{theorem}\cite{D} \label{Thm D} Let J be a homogeneous Gorenstein ideal of codimension three
generated by 2m+1 elements of degrees $d_1\le d_2\leq \cdots \le d_n$ with the resolution (\ref{Jres}). If $c$ is the shift in the last step of the resolution, then \\
\begin{eqnarray*}
c&>&d_i+d_{n-i+2},\  i=2,\ldots,2m+1.\\
\end{eqnarray*}
\end{theorem}

\begin{theorem}\label{thm2.5}
Let $L=(f_1,f_2,f_3,\ldots,f_n)$ be a complete intersection ideal
and $M$ is an ideal minimally generated by $g_1,g_2,\ldots,g_m$ with $n\leq m$. If $L\subset M$, then $degf_i\geq
degg_i$ for all $i=1,2,\cdots,n$.
\begin{proof}
Let $deg f_i=e_i$ and $deg g_j=d_j$ for $i=1,2,\ldots,n$ and
$j=1,2,\ldots,m$.\\
We may assume that $e_1\leq e_2\leq\ldots\leq e_n$, $d_1\leq\
d_2\leq\ldots\leq d_m$. Since $f_1\in L\subset M$,
$f_1=r_{11}g_1+r_{12}g_2+\ldots+r_{1m}g_m$, $r_{1j}\in R$. One of
$g_j$'s is not zero. So $r_{1j}\neq 0$ for some $j$. $f_1$ can be
rewritten as $f_1=r_{1j}g_j+$the other terms. $e_1\geq
deg r_{1j}+d_j$ for some $j$. $d_j\geq d_1$ implies $e_1\geq d_1$.
Since $f_2\in L\subset M$, $f_2=r_{21}g_1+\ldots+r_{2m}g_m$.\
\ \ If $r_{2j}\neq 0$ for some $j\geq 2$, then $f_2=r_{2j}g_j+$the other
terms. $e_2\geq deg r_{2j}+d_j$ for some $j$. we get $e_2\geq d_2$.
We may assume that $r_{2j}=0$ for all $j\geq 2$. Then
$f_2=r_{21}g_1$. So $e_2=deg r_{21}+d_1$. Since $e_1\geq d_1$,
$f_1=r_{11}g_1+r_{1t}g_t+\ldots$ for some $t\neq 1$. If
$f_1=r_{11}g_1$, then $(f_1,f_2)\subset (g_1)$, so codimension of
$(f_1,f_2)=1$ which implies the codimension of $(f_1,\ldots,f_n)$
has at most $n-1$, i.e, codim$(f_1,\ldots,f_n)\leq n-1$. This
contradicts to the codim $I=n$.
If $f_1=r_{11}g_1+r_{1t}g_t+\ldots$ with $r_{1t}\neq 0$ for some $t\neq 1$, then
$e_1\geq d_t$. We may call $t$ by $2$ because $d_2$ is the least
except $d_1$, which implies $e_2\geq e_1\geq d_2$. Inductively, we
get $e_i\geq d_i$ for each $i=1,2,\ldots,n$.
\end{proof}
\end{theorem}

The mutiplicity bounds can be easily established for $R/I$ if the complete intersection $K$ has mutiplicity sufficiently large or sufficiently low and the resolution (\ref{Ires}) is minimal.

\begin{theorem}
If $e(R/K)\leq 3e(R/J)$ then the upperbound in the conjecture 1.1 holds,
\[
e(R/I)\leq \frac{M_1M_2M_3}{6}.
\]
\begin{proof}
We get the maximal shifts in (\ref{Ires}):
\begin{flalign*}
&M_1=max\{e_3,\  e_4\}\\
&M_2=max\{e_2+e_3,\  e_4+d_n\}\\
&M_3=e_1+e_2+e_3-d_1
\end{flalign*}
\begin{flalign*}
6e(R/I)&=6e(R/K)-6e(R/J)\\
&\leq 6e(R/K)-2e(R/K)\\
&=4e(R/K)=4e_1e_2e_3\leq e_3(e_2+e_3)^2\\
&\leq e_3(e_2+e_3)(e_1+e_2+e_3-d_1)\leq M_1M_2M_3
\end{flalign*}
\end{proof}
\end{theorem}

\begin{theorem}
If $e(R/K)\geq 3e(R/J)$ and $e_3 < d_n$ then the lowerbound in the conjecture 1.1 holds,
\[
e(R/I)\geq \frac{m_1m_2m_3}{6}.
\]
\begin{proof}
We get the minimal shifts in (\ref{Ires}):
\begin{flalign*}
&m_1=min\{e_1,\  e_4\}\\
&m_2=min\{e_1+e_2,\  e_4+d_1\}\\
&m_3=e_1+e_2+e_3-d_n
\end{flalign*}
\begin{flalign*}
6e(R/I)&=6e(R/K)-6e(R/J)\\
&\geq 6e(R/K)-2e(R/K)\\
&=4e(R/K)=4e_1e_2e_3\geq e_1(e_1+e_2)^2\\
&\geq e_1(e_1+e_2)(e_1+e_2+e_3-d_n)\geq m_1m_2m_3
\end{flalign*}
\end{proof}
\end{theorem}
\section{The Upper Bound}%

We show that the multiplicity of $R/I$ satisfies the upper bound of
conjecture 1.1. When there is no cancellation at each step in the
resolution (\ref{Ires}), that is exactly the minimal resolution of $R/I$. Furthermore, we have several cases of the minimal free resolutions in which there are cancellations of degrees
between $e_i$'s and $d_i$'s. We can assume $e_1\geq d_1$ or $e_2\geq d_2$ or $e_3\geq d_3$ by theorem ~\ref{thm2.5}. The only cancellations that matter are $e_1=d_1$ or $e_2=d_2$ or $e_3=d_3$. We consider each of these cases separately.

\begin{theorem}
Suppose that $I$ is an almost complete intersection of codim 3 with
a regular sequence $(f_1,f_2,f_3)$ and $J=((f_1,f_2,f_3):I)$ is the
Gorenstein ideal generated in degree $\geq d$ and $M_1$ be the
maximal degree of a minimal generating set of I.\

 Then either
\begin{align*}
M_1&= \sum_{i=1}^{3}deg f_i-c\\
 or\  M_1&= \underset{1\leq i\leq 3}{max}\{deg f_i\}
\end{align*}
  where $c=\frac{1}{m}\sum_{i=1}^{2m+1}d_i$ and $d_i$'s are degrees of
generators for $J$ with $d_1\leq d_2\leq \ldots\leq d_{2m+1}$.
\begin{proof}
Let $e_i=deg f_i$ $1\leq i\leq 4$ and $e_1\leq e_2\leq e_3$. From
the resolution (\ref{Ires}), $M_1=e_1+e_2+e_3-c$ or $M_1=e_3$.
\end{proof}
\end{theorem}

\begin{lemma}
Suppose $I=(K,f_4)$ is an almost complete intersection of codim 3
and $K=(f_1,f_2,f_3)$ is a regular sequence with $degf_i=e_i$ and
$e_1\leq e_2\leq e_3$. Suppose the Gorenstein ideal $J=(K:I)$ is generated by $g_1,\ldots,g_n$ with degree
$g_i=d_i$ and $d_1\leq \cdots \leq d_n$. If $e_1=d_1$, $e_2=d_2$
and $e_3=d_3$, then $g_i$ can be replaced by $f_i$ for each
$i=1,2,3$, i.e, $J=(f_1,f_2,f_3,g_4,\ldots,g_n)$.
\begin{proof}
Since $J=(f_1, f_2, f_3):I$, $f_i\in J$ for $i=1,2,3$. $f_1$ has a linear combination of $g_i's$ of $J$, i.e.,
$f_1=r_1g_1+r_2g_2+\cdots+r_ng_n$ for $r_j\in R$. Then $deg f_1 \geq deg r_1 + deg g_1$. Since $deg f_1 = deg g_1$,
$deg r_1=0$. We can take $r_1=1$. $f_1=g_1+r_2g_2+\cdots+r_ng_n$ enables us to replace $g_1$ by $f_1$ in $J$. So we get $J=(f_1, g_2,\cdots, g_n)$. $f_2\in J$ implies $f_2=s_1f_1 + \sum_{i>1}s_ig_i$, $\tilde{f_2}:= f_2-s_1f_1=\sum_{i>1}s_ig_i$ and $(f_1, \tilde{f_2}, f_3)$ is a regular sequence.  Replace $f_2$ by $\tilde{f_2}$, then $K=(f_1, \tilde{f_2}, f_3)$ and  $\tilde{f_2}=\sum_{i>1}s_ig_i$. Since $deg \tilde{f_2}=deg f_2, deg s_2=0$ and $s_2=1$ can be taken. Now we replace $g_2$ by $\tilde{f_2}$, i.e., $J=(f_1,\tilde{f_2},g_3,g_4,\ldots,g_n)=(f_1,f_2,g_3,g_4,\ldots,g_n)$. With the same arguments, we have $J=(f_1,f_2,f_3,g_4,\ldots,g_n)$.
\end{proof}
\end{lemma}

\begin{lemma}
 \label{Lem1}
If $e_1,e_2,e_3$ are positive integers with $e_1\leq e_2\leq e_3$, then $6e_1e_2e_3\leq
e_1^3+e_1^2e_3+e_2^2e_3+2e_2e_3^2+e_3^3$.

\begin{proof}
Let $e_3=e_2+a$ with $a\geq 0$. Then
\begin{flalign*}
&e_1^3+e_1^2e_2+e_2^2e_3+2e_2e_3^2+e_3-6e_1e_2e_3\\
&=e_1^3+e_1^2(e_2+a)+e_2^2(e_2+a)+2e_2(e_2+a)^2+(e_2+a)^3-6e_1e_2(e_2+a)\\
&=e_1^3+4e_2^2-5e_1^2e_2+(e_1^2+8e_2^2-6e_1e_2)a+5e_2a+a^3\\
&\geq e_1^3+4e_2^3-5e_1^2e_2=(e_2-e_1)(4e_2^2+4e_1e_2-e_1^2)\geq 0 \\
&\ \text{ because } e_1^2+8e_2^2-6e_1e_2\geq 0.
\end{flalign*}
\end{proof}
\end{lemma}
Now we prove the main theorem of this paper. We can
\begin{theorem}
The multiplicity $e(R/I)$ of almost complete intersection $I$ of
codim 3 satisfies the conjectured upper bound.\\
\begin{proof}
\textbf{CASE I}: Suppose that no cancellations occur at each step in (\ref{Ires}). Then
$e_1>d_1$, $e_2>d_2$, $e_3>d_3$. We have the same resolution (2.1) of I as
minimal.

\[
0\,\rightarrow\,
\sum_{i=1}^{2m+1}R(-(\sum_{j=1}^3e_j-d_i))\,\rightarrow\,
\sum_{i=1}^{2m+1}R(-(\sum_{j=1}^3e_j-(c-d_i))) \oplus\sum_{1\le
i<j\le 3}R(-(e_i+e_j)) \
\]
\[
\rightarrow\,R(-(\sum_{j=1}^3e_j-c))\oplus\sum_{j=1}^3R(-e_j)\,\rightarrow\,R\,\rightarrow\,
R/I\rightarrow\, 0.\
\]

Now we get the maximal shifts from the above resolution.
\begin{flalign*}
&M_1=max\{e_3,\  e_4\}\\
&M_2=max\{e_2+e_3,\  e_4+d_n\}\\
&M_3=e_1+e_2+e_3-d_1
\end{flalign*}

where $e_4=e_1+e_2+e_3-c$ by corollary ~\ref{Cor1}.\

We show that $e(R/I)\leq \frac{M_1M_2M_3}{6}$. Let $m_i'$, $1\leq i \leq 3$ denote the minimal shift in the i-th step in the resolution (\ref{Jres}) of $R/J$, then $m_1'=d_1$, $m_2'=c-d_n$, $m_3'=c$ and $6e(R/J)\geq m_1'm_2'm_3'$. Since
$6e(R/I)=6e(R/K)-6e(R/J)\leq 6e_1e_2e_3-m_1'm_2'm_3'$, it
suffices to show that $6e_1e_2e_3-d_1c(c-d_n)\leq M_1M_2M_3$. We verify 4 subcases separately.

\textbf{Subcase 1}: $M_1=e_3$, $M_2=e_2+e_3$\\
 $e_3\geq e_4$ and $e_2+e_3\geq e_4+d_n$ imply $c\geq e_1+e_2$ and
$c-d_n\geq e_1$. We show that $e_3(e_2+e_3)(e_1+e_2+e_3-d_1)\geq
6e_1e_2e_3-d_1c(c-d_n)$.
\begin{flalign*}
&e_3(e_2+e_3)(e_1+e_2+e_3-d_1)-6e_1e_2e_3+d_1c(c-d_n)\\
&\geq e_3(e_2+e_3)(e_1+e_2+e_3-d_1)-6e_1e_2e_3+d_1(e_1+e_2)e_1,\\
&\ \ \text{ because } c\geq e_1+e_2, c-d_n\geq e_1 \\
&=e_3(e_2+e_3)(e_1+e_2+e_3)-d_1(e_2e_3+e_3^2-e_1^2-e_1e_2)-6e_1e_2e_3\\
&\geq e_1e_2e_3+e_2^2e_3+2e_2e_3^2+e_1e_3^2+e_3^3-e_1(e_2e_3+e_3^2-e_1^2-e_1e_2)-6e_1e_2e_3,\\
&\ \ \text{ because } e_1>d_1\\
&=e_3^3+2e_2e_3^2+e_2^2e_3+e_1^2e_2+e_1^3-6e_1e_2e_3\geq 0 \text{ by the lemma}\ \ref{Lem1}.
\end{flalign*}

\textbf{Subcase 2}: $M_1=e_3$, $M_2=e_4+d_n$\\
 $e_3\geq e_4$ and $e_4+d_n\geq e_2+e_3$ imply $c\geq e_1+e_2$ and
$c-d_n\leq e_1$.We show that $e_3(e_1+e_2+e_3-c+d_n)(e_1+e_2+e_3-d_1)\geq 6e_1e_2e_3-d_1c(c-d_n)$.
\begin{flalign*}
&e_3(e_1+e_2+e_3-c+d_n)(e_1+e_2+e_3-d_1)-6e_1e_2e_3+d_1c(c-d_n)\\
&=e_3(e_1+e_2+e_3)(e_1+e_2+e_3-d_1)-(c-d_n)\{e_3(e_1+e_2+e_3-d_1)-d_1c\}\\
&\ -6e_1e_2e_3\\
&\geq e_3(e_1+e_2+e_3)(e_1+e_2+e_3-d_1)-e_1(e_3-d_1)(e_1+e_2+e_3)-6e_1e_2e_3\\
&\ \ \text{ because } c\geq e_1+e_2,\  c-d_n\leq e_1 \\
&=e_3(e_1+e_2+e_3)^2-e_1e_3(e_1+e_2+e_3)-d_1(e_3-e_1)(e_1+e_2+e_3)-6e_1e_2e_3\\
&\geq e_3(e_2+e_3)(e_1+e_2+e_3)-e_1(e_3-e_1)(e_1+e_2+e_3)-6e_1e_2e_3 \text{ since }e_1>d_1\\
&= e_3^3+2e_2e_3^2+e_2^2e_3+e_1^2e_2+e_1^3-6e_1e_2e_3\geq 0 \text{ by the lemma}\ \ref{Lem1}.
\end{flalign*}

\textbf{Subcase 3}: $M_1=e_4$, $M_2=e_2+e_3$\\
$e_4\geq e_3$ and $e_2+e_3\geq e_4+d_n$ imply $c\leq e_1+e_2$ and
$c-d_n\geq e_1.$\\
We show that $(e_1+e_2+e_3-c)(e_2+e_3)(e_1+e_2+e_3-d_1)\geq
6e_1e_2e_3-d_1c(c-d_n)$.
\begin{flalign*}
&(e_1+e_2+e_3-c)(e_2+e_3)(e_1+e_2+e_3-d_1)+d_1c(c-d_n)-6e_1e_2e_3\\
&=(e_2+e_3)(e_1+e_2+e_3)(e_1+e_2+e_3-c-d_1)+d_1c(e_2+e_3+c-d_n)-6e_1e_2e_3\\
&\geq (e_2+e_3)(e_1+e_2+e_3)(e_3-d_1)+d_1c(e_1+e_2+e_3)-6e_1e_2e_3\\
&\ \ \text{ because } e_1+e_3\geq c,\  c-d_n\geq e_1\\
&= e_3(e_2+e_3)(e_1+e_2+e_3)-d_1\{(e_2+e_3)(e_1+e_2+e_3)-c(e_1+e_2+e_3)\}-6e_1e_2e_3\\
&\geq e_3(e_2+e_3)(e_1+e_2+e_3)-e_1(e_3-e_1)(e_1+e_2+e_3)-6e_1e_2e_3\\
& \ \ \text{ because } e_1>d_1,\  c\leq e_1+e_2\\
&= e_3^3+2e_2e_3^2+e_2^2e_3+e_1^2e_2+e_1^3-6e_1e_2e_3\geq 0 \text{
by the lemma}\ \ref{Lem1}.
\end{flalign*}

\textbf{Subcase 4}: $M_1=e_4$, $M_2=e_4+d_n$\\
$e_4\geq e_3$ and $e_4+d_n\geq e_2+e_3$ imply $c\leq e_1+e_2$ and
$c-d_n\leq e_1.$\\
We show that
$(e_1+e_2+e_3-c)(e_1+e_2+e_3-c+d_n)(e_1+e_2+e_3-d_1)\geq 6e_1e_2e_3-d_1c(c-d_n)$.
\begin{flalign*}
&(e_1+e_2+e_3-c)(e_1+e_2+e_3-c+d_n)(e_1+e_2+e_3-d_1)+d_1c(c-d_n)-6e_1e_2e_3\\
&=(e_1+e_2+e_3-c)(e_1+e_2+e_3-d_1)(e_1+e_2+e_3)\\
&\ \ -(c-d_n)\{(e_1+e_2+e_3-c)(e_1+e_2+e_3-d_1)-d_1c\}-6e_1e_2e_3\\
&\geq e_3(e_2+e_3)(e_1+e_2+e_3)-e_1(e_3-e_1)(e_1+e_2+e_3)-6e_1e_2e_3,\\
&\ \ \text{ because } e_1+e_2\geq c,\  c-d_n\leq e_1\\
&= e_3^3+2e_2e_3^2+e_2^2e_3+e_1^2e_2+e_1^3-6e_1e_2e_3\geq 0 \text{
by the lemma} ~\ref{Lem1}.
\end{flalign*}

\textbf{CASE II}: Suppose that there is the cancellation of $e_1=d_1$ in the
resolution (\ref{Jres}) of $I$. Then, the minimal resolution is
the following.

\[
0\,\rightarrow\,
\sum_{i=2}^{n}R(-(e_1+e_2+e_3-d_i))\,\rightarrow\,
\sum_{i=1}^{n}R(-(e_4-d_i))\oplus R(-(e_1+e_2))\oplus R(-(e_1+e_3))
\]

\[
\rightarrow\,\sum_{i=1}^{4}R(-(e_i))\,\rightarrow\,R\,\rightarrow\,
R/I\rightarrow\, 0\
\]

Here, $e_2$ must be strictly greater than $d_2$. Otherwise, there is
more cancellation of $e_2=d_2$ for which we have to consider later.\\
Then we have the maximal shifts:
\begin{flalign*}
&M_1=max\{e_3,\  e_4\}\\
&M_2=max\{e_1+e_3,\ e_4+d_n\}\\
&M_3=e_1+e_2+e_3-d_2
\end{flalign*}

If $M_1=e_4$, then $e_1+e_2+e_3-c>e_3$. We get $e_1+e_2>c$, which
always results in $M_2=e_4+d_n$ because $e_4+d_n-(e_1+e_3)=e_1+e_2+e_3-c+d_n-e_1-e_3=e_2$
+$d_n-c>e_2+e_1-c>0.$ We just have 3 cases of $(M_1=e_3,\textbf{ } M_2=e_1+e_3)$, $(M_1=e_3,\textbf{ }
M_2=e_4+d_n)$ and $(M_1=e_4,\textbf{ } M_2=e_4+d_n)$. We show the upper bound
of the conjecture (1.1) with when $M_2=e_1+e_3$ or $M_2=e_4+d_n$.\\

\textbf{Subcase 1}: $M_2=e_1+e_3$\\
$M_2=e_1+e_3\geq e_4+d_n$ implies $c-d_n\geq e_2$\\
In \cite{M-N-R}, Migliore, Nagel and R\"{o}mer achieved the stronger bound of
multiplicity of a Gorenstein ideal of codim 3 as following\\
\begin{flalign*}
6e(R/J)\geq e_1c(c-d_n)+2e_1^2(d_n-e_1)
\end{flalign*}

Hence $6e(R/I)=6e_1e_2e_3-6e(R/J)\leq
6e_1e_2e_3-e_1c(c-d_n)-2e_1^2(d_n-e_1)$. Since $M_1\geq e_3$,
\begin{flalign*}
&M_1M_2M_3-6e(R/I)\\
&\geq e_3(e_1+e_3)(e_1+e_2+e_3-d_2)-6e_1e_2e_3+e_1c(c-d_n)+2e_1^2(d_n-e_1)\\
&\geq e_3(e_1+e_3)^2+e_1(e_2+d_n)e_2+2e_1^2(d_n-e_1)-6e_1e_2e_3.
\end{flalign*}

We verify that the right above expression is nonnegative.\

If $d_n\geq e_2$, then the above expression is greater than
$e_3(e_1+e_3)^2+2e_1e_2^2+2e_1^2e_2-2e_1^3-6e_1e_2e_3$.\

Let $e_3=e_2+a$, $a\geq 0$ then the above can be rewritten as
\begin{flalign*}
&e_3(e_1+e_3)^2+2e_1e_2^2+2e_1^2e_2-2e_1^3-6e_1e_2e_3\\
&=3e_1^2e_2+e_2^3-2e_1^3-e_1e_2^2+(e_1^2+3e_2^2-2e_1e_2)a+(2e_1+3e_2)a^2+a^3\\
&\geq
3e_1^2e_2+e_2^3-2e_1^3-e_1e_2^2=(e_2^3-e_1^3)-e_1(e_1^2-3e_1e_2+2e_2^2)\\
&=(e_2-e_1)(e_2^2+e_1e_2+e_1^2-e_1(2e_2-e_1))=(e_2-e_1)(e_2^2-e_1e_2+2e_1^2)\geq
0
\end{flalign*}

If $d_n<e_2$, then
\begin{flalign*}
&e_3(e_1+e_3)(e_1+e_2+e_3-d_2)+e_1c(c-d_n)+2e_1^2(d_n-e_1)-6e_1e_2e_3\\
&\geq
e_3(e_1+e_3)(e_1+e_2+e_3-d_2)+e_1e_2(e_2+d_n)+2e_1^2(d_n-e_1)-6e_1e_2e_3\\
&\geq
e_3(e_1+e_3)(e_1+e_2+e_3-d_n)+e_1e_2(e_2+d_n)+e_1^2(d_n-e_1)-6e_1e_2e_3\\
&=e_3(e_1+e_3)(e_1+e_2+e_3)+e_1e_2^2-e_1^3-6e_1e_2e_3-d_n(e_1e_3+e_3^2-e_1e_2-e_1^2)\\
&\geq
e_3(e_1+e_3)(e_1+e_2+e_3)+e_1e_2^2-e_1^3-6e_1e_2e_3-e_2(e_1e_3+e_3^2-e_1e_2-e_1^2)\\
&=e_3^3+2e_1e_3^2+e_1^2e_3+2e_1e_2^2+e_1^2e_2-e_1^3-6e_1e_2e_3
\end{flalign*}
Now let's take $e_3$ by $e_2+a$, $a\geq 0$, then the right above
expression is
\begin{flalign*}
&(e_2+a)^3+2e_1(e_2+a)^2+e_1^2(e_2+a)+2e_1e_2^2+e_1^2e_2-e_1^3-6e_1e_2(e_2+a)\\
&=e_2^3-e_1^3-2e_1e_2^2+2e_1^2e_2+a^3+(3e_2+2e_1)a^2+(3e_2^2-2e_1e_2+e_1^2)a\\
&\geq(e_2-e_1)\{(e_2-e_1)^2+e_1e_2\}\geq 0
\end{flalign*}

\textbf{Subcase 2}: $M_2=e_4+d_n$\\
$e_4+d_n\geq e_1+e_3$ implies $e_2\geq c-d_n$.\\
Since $M_1\geq e_3$, we show that $e_3(e_1+e_2+e_3-c+d_n)(e_1+e_2+e_3-d_2)\geq
6e_1e_2e_3-e_1c(c-d_n)-2e_1^2(d_n-e_1)$.\\ \\
If $d_n\geq e_2$, then $c>d_2+d_n\geq d_2+e_2\geq e_1+e_2$ and
\begin{flalign*}
&e_3(e_1+e_2+e_3-c+d_n)(e_1+e_2+e_3-d_2)+e_1c(c-d_n)+2e_1^2(d_n-e_1)-6e_1e_2e_3\\
&\geq
e_3(e_1+e_3)(e_1+e_2+e_3-c+d_n)+e_1(e_1+e_2)(c-d_n)+2e_1^2(e_2-e_1)-6e_1e_2e_3\\
&\geq e_3(e_1+e_3)(e_1+e_2+e_3)-(c-d_n)\{e_3(e_1+e_3)-e_1(e_1+e_2)\}+2e_1^2(e_2-e_1)\\
&\ \ -6e_1e_2e_3\\
&\geq e_3(e_1+e_3)(e_1+e_2+e_3)-e_2\{e_3(e_1+e_3)-e_1(e_1+e_2)\}+2e_1^2(e_2-e_1)-6e_1e_2e_3\\
&\geq e_3(e_1+e_3)^2+e_1e_2(e_1+e_2)+2e_1^2(e_2-e_1)-6e_1e_2e_3
\end{flalign*}
Let $e_3=e_2+a$, then $a\geq 0$. Replace $e_3$ by $e_2+a$ to the right above then we have
\begin{flalign*}
&(e_2+a)(e_1+e_2+a)^2+3e_1^2e_2+e_1e_2^2-2e_1^3-6e_1e_2(e_2+a)\\
&=e_2(e_1+e_2)^2+3e_1^2e_2+e_1e_2^2-2e_1^3-6e_1e_2^2+\{(e_1+e_2)^2+2e_2(e_1+e_2)-6e_1e_2\}a\\
&\ +(2e_1+3e_2)a^2+a^3\geq (e_2-e_1)^3\geq 0\\
&\ \ \text{ because }(e_1+e_2)^2+2e_1(e_1+e_2)-6e_1e_2=(e_2-e_1)(3e_2-e_1)\geq 0 \text{ and } a\geq 0.
\end{flalign*}
If $d_n<e_2$, then $2e_2>e_2+d_n>c>d_2+d_n$. The last
inequality results from theorem \ref{Thm D}.
\begin{flalign*}
&e_3(e_1+e_2+e_3-c+d_n)(e_1+e_2+e_3-d_2)+e_1c(c-d_n)+2e_1^2(d_n-e_1)-6e_1e_2e_3\\
&=e_3(e_1+e_2+e_3)(e_1+e_2+e_3-d_2)\\
&\ \ \ -(c-d_n)\{e_3(e_1+e_2+e_3-d_2)-e_1c\}+2e_1^2(d_n-e_1)-6e_1e_2e_3\\
&\geq e_3(e_1+e_2+e_3)(e_1+e_2+e_3-d_2)\\
&\ \ \ -e_2\{e_3(e_1+e_2+e_3-d_2)-e_1(d_2+d_n)\}+2e_1^2(d_n-e_1)-6e_1e_2e_3,\\
&\ \ \text{ because }
e_3(e_1+e_2+e_3-d_2)-e_1(d_2+d_n)>e_3(e_1+e_2+e_3-d_2)-e_1c\\
&\ \ \ >e_3(e_1+e_3)-2e_1e_2>0.\\
&\geq
e_3(e_1+e_3)(e_1+e_2+e_3-d_2)+e_1e_2(d_2+d_n)+2e_1^2(d_n-e_1)-6e_1e_2e_3\\
&=e_3(e_1+e_3)(e_1+e_2+e_3)-d_2\{e_3(e_1+e_3)-e_1e_2\}+e_1e_2d_n+2e_1^2(d_n-e_1)-6e_1e_2e_3\\
&\geq
e_3(e_1+e_3)(e_1+e_2+e_3)-d_n\{e_3(e_1+e_3)-e_1e_2\}+e_1e_2d_n+2e_1^2(d_n-e_1)-6e_1e_2e_3\\
&=e_3(e_1+e_3)(e_1+e_2+e_3-d_n)+2e_1e_2d_n+2e_1^2(d_n-e_1)-6e_1e_2e_3.
\end{flalign*}

We claim that the right above expression is greater than
$e_3(e_1+e_3)^2+e_1e_2(e_1+e_2)+e_1^2(e_2-e_1)-6e_1e_2e_3$ which is
already proved to be nonnegative in the end of case $d_n\geq e_2$.\\
To verify our claim, it suffices to show
$e_3(e_1+e_3)(e_2-d_n)+e_1e_2(2d_n-e_1-e_2)+e_1^2(2d_n-2e_1-e_2+e_1)\geq
0$. Let $e_2=d_n+x$, then $0<x\leq e_3-d_n$.

\begin{flalign*}
&e_3(e_1+e_3)(e_2-d_n)+e_1e_2(2d_n-e_1-e_2)+e_1^2(2d_n-e_1-e_2)\\
&=e_3(e_1+e_3)x+e_1(d_n+x)(2d_n-e_1-d_n-x)+e_1^2(2d_n-e_1-d_n-x)\\
&=e_1(d_n+e_1)(d_n-e_1)+x(e_3^2+e_1e_3-2e_1^2-e_1x)\\
&\geq e_1(d_n+e_1)(d_n-e_1)+x\{(e_3^2+e_1e_3-2e_1^2-e_1(e_3-d_n)\}\\
&=e_1(d_n+e_1)(d_n-e_1)+x(e_3^2-2e_1^2+e_1d_n)\geq 0.\\
\end{flalign*}

\textbf{CASE III }: Suppose that $e_1=d_1$, $e_2=d_2$, $e_3>d_3$.\
Denote $d_1$ by $e_1$ and $d_2$ by $e_2$.
\[
0\,\rightarrow\,
\sum_{i=3}^{n}R(-(e_1+e_2+e_3-d_i))\,\rightarrow\,
\sum_{i=1}^{n}R(-(e_4+d_i))\oplus R(-(e_1+e_2))
\]
\[
\rightarrow\,\sum_{i=1}^{4}R(-(e_i))\,\rightarrow\,R\,\rightarrow\,
R/I\rightarrow\, 0\
\]
$e_3-e_4=e_3-(e_1+e_2+e_3-c)=c-(e_1+e_2)>c-(e_2+d_n)>0$ applied to theorem~\ref{Thm D}. This gives us the following maximal shifts.
\begin{flalign*}
M_1 &= e_3\\
M_2 &= max\{e_1+e_2,\  e_4+d_n\}\\
M_3 &= e_1+e_2+e_3-d_3
\end{flalign*}

\textbf{Subcase 1}: $M_2=e_1+e_2$\\
$e_1+e_2\geq e_4+d_n$ implies $e_3\leq c-d_n$.
\begin{flalign*}
&e_3(e_1+e_2)(e_1+e_2+e_3-d_3)+e_1c(c-d_n)-6e_1e_2e_3\\
&\geq e_3(e_1+e_2)^2+e_1(e_3+d_n)e_3-6e_1e_2e_3\\
&\geq e_3(e_1+e_2)^2+e_1(e_2+e_3)e_3-6e_1e_2e_3,\text{ because } d_n\geq e_2\\
&=e_1^2e_3+e_2^2e_3+e_1e_3^2-3e_1e_2e_3=e_3(e_1^2+e_2^2+e_1e_3-3e_1e_2)\\
&=e_3\{(e_1-e_2)^2+e_1e_3-e_1e_2\}=e_3\{(e_1-e_2)^2+e_1(e_3-e_2)\}\geq 0
\end{flalign*}

\textbf{Subcase 2}: $M_2=e_4+d_n$\\
$M_2=e_4+d_n>e_1+e_2$ implies $e_3>c-d_n$.We know that $c>e_2+d_n$, so $e_3\geq c-d_n>e_2$.\\
If $e_3\geq d_n$, then we show that
$e_3(e_1+e_2+e_3-c+d_n)(e_1+e_2+e_3-d_3)\geq 6e_1e_2e_3-e_1c(c-d_n)$.
\begin{flalign*}
&e_3(e_1+e_2+e_3-c+d_n)(e_1+e_2+e_3-d_3)+e_1c(c-d_n)-6e_1e_2e_3\\
&=e_3(e_1+e_2+e_3)(e_1+e_2+e_3-d_3)\\
&\ \ -(c-d_n)\{e_3(e_1+e_2+e_3-d_3)-e_1c\}-6e_1e_2e_3\\
&\geq e_3(e_1+e_2+e_3)(e_1+e_2+e_3-d_3)-e_3\{e_3(e_1+e_2+e_3-d_3)-e_1c\}-6e_1e_2e_3\\
&\ \text{ because } e_3(e_1+e_2+e_3-d_3)-e_1c\geq
e_3(e_1+e_2)-e_1(e_3+d_n)\geq e_2e_3-e_1d_n>0\\
&=e_3(e_1+e_2)(e_1+e_2+e_3-d_3)+e_1e_3c-6e_1e_2e_3\\
&>e_3(e_1+e_2)^2+e_1e_3(2e_2)-6e_1e_2e_3 \text{  because  }
c>e_2+d_n>2e_2\\
&=e_3(e_2-e_1)^2\geq 0
\end{flalign*}

If $e_3<d_n$, then $e_2+e_3<e_2+d_n<c\leq e_3+d_n<2d_n$.\\
We show that
$e_3(e_1+e_2+e_3-c+d_n)(e_1+e_2+e_3-d_3)+e_1c(c-d_n)+2e_1^2(d_n-e_1)\geq
6e_1e_2e_3$.
\begin{flalign*}
&e_3(e_1+e_2+e_3-c+d_n)(e_1+e_2+e_3-d_3)+e_1c(c-d_n)+2e_1^2(d_n-e_1)-6e_1e_2e_3\\
&\geq e_3(e_1+e_2)^2+e_1(e_2+d_n)e_2+2e_1^2(d_n-e_1)-6e_1e_2e_3\\
&\geq e_3(e_1+e_2)^2+e_1e_2(e_2+e_3)+2e_1^2(e_3-e_1)-6e_1e_2e_3\\
&=3e_1^2e_3+e_2e_3^2+e_1e_2^2-2e_1^3-3e_1e_2e_3\\
\end{flalign*}

Let $e_3=e_2+a$, $a>0$ then the right above expression is
\begin{flalign*}
&3e_1^2(e_2+a)+e_2(e_2+a)^2+e_1e_2^2-2e_1^3-3e_1e_2(e_2+a)\\
&=3e_1^2e_2+e_2^3-2e_1^3-2e_1e_2^2+(3e_1^2+2e_2^2-3e_1e_2)a+e_2a^2\\
&\geq
(e_2-e_1)(e_2^2-e_1e_2+2e_1^2)\geq 0
\end{flalign*}

\textbf{CASE IV}: Suppose the case of $e_1=d_1$, $e_2=d_2$, $e_3=d_3$.
By theorem \ref{thm2.5}  we can take $f_i = g_i, 1\leq i \leq 3$.  Thus 
$J = ( (f_1, f_2, f_3): f_4)$ and $I = ( (f_1, f_2, f_3): J)$.  Since J is a homogeneous Gorenstein ideal of  height three,  there exists a skew symmetric matrix $\phi $ of size $n= 2m+1$ such that  $ J $ is the ideal of $2m\times 2m$ order pfaffians of $\phi $.  The degree matrix of $\phi $ is the following. 

$$\left[
\begin{matrix}
2r_1 & r_1+r_2 & \cdots & r_1+r_n \\
r_1+r_2 & 2r_2 & \cdots & r_2+r_n \\
\vdots & \vdots &       & \vdots \\
r_1+r_n & r_2+r_n & \cdots & 2r_n
\end{matrix}
\right]$$

 Let $r=r_1+r_2+\cdots +r_n$, Then homogeneous generators, $g_i$  of $J$ are of degrees, $deg g_i = r-r_i$.  We may assume that $r_1\geq r_2\geq \cdots \geq r_n$. The determinant of $\phi$ is a homogeneous polynomial of degree $c=2(r_1+r_2+\cdots +r_n)$.   We have $d_i=r-r_i$ for all $i=1,2,\cdots ,n$ and $e_4=d_1+d_2+d_3-c=r-r_1-r_2-r_3$. Since $e_4$ and $c$ is a positive integer, $r_1,r_2,r_3,r_4$ and $r$ are positive.  With this notation, the free resolution of $R/I$  is as follows:
 
\[
0\,\rightarrow\,
\sum_{i=4}^{n}R(-(2r-r_1-r_2-r_3+r_i))\,\rightarrow\,
\sum_{i=1}^{n}R(-(2r-r_1-r_2-r_3-r_i))
\]

\[
\rightarrow\,\sum_{i=1}^{3}R(-(r-r_i))\oplus R(-\sum_{j=4}^{n}r_j)\,\rightarrow\,R\,\rightarrow\,
R/I\rightarrow\, 0\
\]

Let $T = 2r-r_1-r_2-r_3$.  We have the maximal and minimal shifts.
\begin{alignat*}{3}
        M_1&=r-r_3& \qquad   m_1&=T-r \\
            M_2&=T-r_n& \qquad   m_2&=T-r_1\\
            M_3&=T+r_4& \qquad   m_3&=T+r_n
\end{alignat*}

 Using this resolution to compute the multiplicity of $R/I$, we get 
 
 \begin{flalign*}
6e(R/I)&=\sum_{i=4}^{n}(T-r_i)^3-\sum_{i=1}^{n}(T-r_i)^3+\sum_{i=1}^3(r-r_i)^3+(T-r)^3\\
     &= 3T(\sum_{i=4}^{n}r_i^2-\sum_{i=1}^{n}r_i^2)+\sum_{i=4}^{n}r_i^3+\sum_{i=1}^{n}r_i^3+\sum_{i=1}^{3}(r-r_i)^3+(T-r)^3\\
     &=-3T\sum_{i=1}^{3}r_i^2+2\sum_{i=4}^{n}r_i^3+3r^2(T-r)+3r\sum_{i=1}^{3}r_i^2+(T-r)^3\\
     &=(T-r)(3r^2 - 3\sum_{i=1}^{3}r_i^2)+2\sum_{i=4}^{n}r_i^3+(T-r)^3
\end{flalign*}
Let $ \Delta :=6e(R/I)-M_1M_2M_3$.  We will show that  $\Delta \leq 0$.
\begin{flalign*}
   \Delta &=(T-r)(3r^2 - 3\sum_{i=1}^{3}r_i^2)+2\sum_{i=4}^{n}r_i^3+(T-r)^3-(r-r_3)(T-r_n)(T+r_4)\\
&=(T-r)(3r^2 - 3\sum_{i=1}^{3}r_i^2)+2\sum_{i=4}^{n}r_i^3+(T-r)^3-(T-r)(T-r_n)(T+r_4)\\
&\ \ \ -(r_1+r_2){(r-r_n) +(T-r)}(T+r_4), \textbf{ }\textbf{ } \text{ since } r-r_3=T-r+r_1+r_2 \\
&=(T-r)(3r^2 - 3\sum_{i=1}^{3}r_i^2)+2\sum_{i=4}^{n}r_i^3+(T-r)^3-(T-r)(T-r_n)(T+r_4)\\
&\ \ \ -(r_1+r_2)(T-r)(T+r_4)-(r-r_n)(r_1+r_2)T-(r-r_n)(r_1+r_2)r_4
 \end{flalign*}
 
   But  $2\sum_{i=4}^{n}r_i^3-(r_1+r_2)r_4(T-r) \leq \sum_{i=4}^{n}(2r_i)r_i r_i-(r_1+r_2)r_4(T-r)\leq 0.$  We get, 

\begin{flalign*}\Delta&\leq(T-r)(3r^2 - 3\sum_{i=1}^{3}r_i^2)+(T-r)^3\\
&\ \ \ -(T-r)(T-r_n)(T+r_4)-(r_1+r_2)(T-r)T-(r-r_n)(r_1+r_2)T\\
&\ \ \ -(r-r_n)(r_1+r_2)r_4 \\
&=(T-r)(3r^2-3\sum_{i=1}^{3}r_i^2)+(T-r)^3-(T-r)(T-r_n)(T+r_4)\\
&\ \ \ -(r_1+r_2)(T-r_n)T-(r-r_n)(r_1+r_2)r_4 \\
&=(T-r)\left\{2r^2-3\sum_{i=1}^{3}r_i^3-2r(T-r)-T(r_4-r_n)+r_4r_n\right\}\\
&\ \ \ -(r_1+r_2)(T-r_n)T-(r-r_n)(r_1+r_2)r_4\\
&=2r^2(T-r)-2r(T-r)^2-(T-r)\left\{T(r_4-r_n)+ 3\sum_{i=1}^{3}r_i^3-r_4r_n\right\}\\
&\ \ \ -(r_1+r_2)(T-r_n)T-(r-r_n)(r_1+r_2)r_4\\
&=2r(T-r)(r_1+r_2+r_3)-(T-r)\left\{T(r_4-r_n)+ 3\sum_{i=1}^{3}r_i^3-r_4r_n\right\}\\
&\ \ \ -(r_1+r_2)(T-r_n)T-(r-r_n)(r_1+r_2)r_4.
\end{flalign*}

$T-r=r-(r_1+r_2+r_3)\leq r-r_n$ implies $2(T-r)\leq T-r+r-r_n=T-r_n$.   Now, 
$T-r = r_4+\cdots +r_n \geq r_n$ whether $r_n$  is positive or negative since $T-r$ is 
positive.   Now, 
\begin{flalign*}
&2r(T-r)(r_1+r_2+r_3) \\
&=2r(T-r)(r_1+r_2)+2r r_3(T-r)\leq r(r_1+r_2)(T-r_n)+2rr_3(T-r)\\
&\leq r(r_1+r_2)(T-r_n)+r(r_1+r_2)(T-r) \text{ because } 2r_3 < r_1+r_2\\
&\leq r(r_1+r_2)(T-r_n)+(r_1+r_2)(T-r_n)(T-r)= (r_1+r_2)(t-r_n)T
\end{flalign*}

Thus 

\begin{flalign*}
\Delta&\leq -(T-r)\left\{T(r_4-r_n)+ 3\sum_{i=1}^{3}r_i^3-r_4r_n\right\}-(r-r_n)(r_1+r_2)r_4\\
&\leq -(T-r)\left\{T(r_4-r_n)+ 2r_1^2+3r_2^2+3r_3^2\right\}-(r-r_n)(r_1+r_2)r_4\leq 0.
\end{flalign*}
\end{proof}
\end{theorem}


\begin{remark} When $e_1=d_1, e_2=d_2, e_3= d_3$, the structure of the ideal I and 
its resolution can be completely determined.  There exists a homogeneous 
skew symmetric  matrix $\phi   = (x_{ij})$ giving the Gorenstein ideal $J$, with $\tilde{\phi}$ denoting $\phi$ with the top three rows deleted, the resolution of $R/I$ is

\begin{eqnarray*} 
0\,\rightarrow\,
R^{n-3} \ \stackrel{\tilde{\phi}} \rightarrow\   R^n \ \stackrel {\psi }
\rightarrow\,   R^4 \,\stackrel {(f_4,
f_1,
f_2,
f_3)^t} \longrightarrow\,   R\,\rightarrow\,R/I\,\rightarrow\, 0.
\end{eqnarray*}

where $\psi $ is the $n \times 4$ matrix
$$\left[
\begin{matrix}
f_1     &         &                   &  \\
f_2     &         &-f_4 I_{3\times 3} &  \\
f_3     &         &                   &  \\
g_4     & t_{14}  & t_{24}            & t_{34}\\
\vdots  & \vdots  &  \vdots           &  \vdots\\ 
g_n     & t_{1n}  &  t_{2n}           & t_{3n}
\end{matrix}
\right]
$$
\end{remark}

\begin{example} $I=(x^7, y^8+z^8, x^3y^6+x^5z^4+yz^8, y^3z^3)$ is an
almost complete intersection of codim 3.
 Let $f_1=x^7, f_2=y^8+z^8,
f_3=x^3y^6+x^5z^4+yz^8$ and $f_4=y^3z^3$, then $K=(f_1, f_2, f_3)$ is
a regular sequence.\

The Gorenstein ideal $ J=((f_1, f_2, f_3):I)$ is $(x^3y^3z-y^6z+x^2z^5,
x^6z+xyz^5, x^7, x^5y^2z-x^3z^5+y^3z^5, y^8+z^8, x^5y^3-x^2y^6,
x^3y^6+x^5z^4-yz^8)$ of $7$ generators and $c=\frac{1}{3}(7+7+7+
8+8+8+9)=18$. $e(R/I)=270=504-234=e(R/K)-e(R/J)$ in lemma ~\ref{Lem2} and $deg f_4=6=deg f_1+deg f_2+deg f_3-c$ in corollary ~\ref{Cor1}. Since $deg f_1=deg g_1=7$, this example is for the case II in the previous theorem. The maximal and minimal shifts for $R/I$ are $M_1=e_3=9, M_2=e_1+e_3=16, M_3=e_1+e_2+e_3-d_2=17$ and $m_1=e_4=6, m_2=e_4+d_1=13, m_3=e_1+e_2+e_3-d_n=15$. $1170=m_1m_2m_3 \leq 6e(R/I) \leq M_1M_2M_3=2448$.
\end{example}

\section{Ideal Linked to Complete Intersection}%

In this section we use some techniques of Gold-Schenck-Srinivasan \cite{G-S-S} to give bounds for the almost complete intersection ideal linked to a complete intersection in one step and give some partial results towards the proving the multiplicity conjecture for these ideals. 
I is a complete intersection generated by homogeneours elements of 
degrees $d_1\le d_2\le \cdots \le d_n$.  K is  a regular sequence of 
length n  formed by homogeneous elements  of degrees $e_1\le 
e_2\le \cdots \le e_n$.  By theorem \ref{thm2.5} without loss of generality, we can take $e_i \ge 
d_i$  since $K \subset I$.   

Let $J= (K:I)$.  Let $\bold{F}$ be the minimal resolution of $R/K$ 
and $\bold{K}$ be the minimal resolution of $R/I$. Let $\phi : \textbf{F} \rightarrow \textbf{K}$ be induced by $K \subset I$, then the dual of the mapping cone of $\phi$, $M(\phi^*)$ is a resolution of $R/J$. This will be minimal if there are no cancellations. We consider this situation where $M(\phi^*)$ is the minimal resolution of $R/I$.

 \textbf{F} is a Koszul complex with degrees of the generators  $d_1\le d_2\le \cdots \le d_n$ and \textbf{K} is a Koszul complex with degrees of the generators  $e_1\le e_2\le \cdots \le e_n$ and let $\alpha = \sum _i e_i$. Then the maximal and minimal shifts in $\textbf{P}=M(\phi^*)$ are 
 
\begin{flalign*}
&M_i = max\{\sum _{t \geq n-i+1}e_t,\  \alpha - \sum_{t \leq n-i+1} d_t\},\  1\leq i\leq n-1\\
&M_n = \alpha -d_1.\\
&m_i = min\{\sum _{t \leq i}e_t,\  \alpha - \sum_{t \geq i} d_t\},\  1\leq i\leq n-1\\
&m_n = \alpha -d_n.
\end{flalign*}

The multiplicity of $R/J$ is
\begin{eqnarray*}
e(R/J)=\prod_{i=1}^ne_i-\prod_{i=1}^nd_i=e(R/K)-e(R/I).\\ 
\end{eqnarray*}

We show with some conditions that
\begin{equation*}
\frac{\prod_{i=1}^nm_i}{n!} \leq e(R/J) \leq \frac{\prod_{i=1}^nM_i}{n!} .\\
\end{equation*}

\begin{theorem}
The upper bound of the conjecture 1.1 holds if $\sum_{i=2}^n(e_i-e_1)\geq
d_1.$
\begin{proof}
 Since $M_i = max\{ \sum _{t \geq n-i+1}e_t,\  \alpha - \sum_{t \leq n-i+1}d_t\}$, $M_i \geq ie_{n-i+1}$, for $1\leq i \leq n-1$. 
 
\begin{flalign*}
 &\prod_{i=1}^n M_i - n!e(R/J)\\ 
 &=\prod_{i=1}^n M_i - n!\prod_{i=1}^n e_i +n!\prod_{i=1}^n d_i\\
 &\ge  M_n \prod _{i<n} ie_{n - i+1} - ne_1 ( n-1)! \prod_{i>1} e_i +n!\prod_{i=1}^n d_i \\
 &\ge (\alpha - d_1 - ne_1) (n-1)!\prod _{i>1} e_i +n! \prod_{i=1}^n d_i. 
\end{flalign*}

This is non negative if $\alpha -d_1 \ge ne_1$.  $\alpha-ne_1 \ge d_1$ if $\sum _i (e_i-e_1) \ge d_1$.
\end{proof}
\end{theorem}

Now we look at the case where the regular sequence linking the two ideals is generated in a single degree $e$, so that  $e=e_1=e_n$.   Then $e \geq d_i$ for all $i$ by the theorem ~\ref{thm2.5} and 

\begin{flalign*}
&M_i = max\{ie,\  ne -\sum_{j=1}^{n-i+1}d_j\}, 1\leq i\leq n-1\\ 
&M_n = ne-d_1\\
&m_i = min\{ie,\  ne -\sum_{j=i}^{n}d_j\}, 1\leq i\leq n-1\\
&m_n = ne-d_n.\\
\end{flalign*}

\begin{lemma}\label{lemm1}
  If $m_k=ke$ for some $k$, then $m_i=ie$ for all $i\leq k$.
\begin{proof}
Since $m_k=ke$, $ne-(d_k+\cdots +d_n)\geq ke$.
$ne-(d_{k-1}+d_k+\cdots +d_n)=ne-(d_k+\cdots +d_n)-d_{k-1}\geq
ke-d_{k-1}\geq ke-e=(k-1)e$  because $d_{k-1}\leq e$, so
$m_{k-1}=(k-1)e$. Repeat this process, $m_i=ie$ $i=1,2,\cdots ,k$.
\end{proof}
\end{lemma}

\begin{theorem}
If $m_{n-1}=(n-1)e$, then the lower bound holds.
\begin{proof}
Since $m_{n-1} = (n-1)e \leq ne-d_{n-1} -d_n$, We get $e \geq d_n + d_{n-1}$.
By the lemma\ \ref{lemm1}, $m_i=ie$ for all $i\leq k$.\\
\begin{eqnarray*}
 e&\geq& \frac{d_1+\cdots +d_n}{n-1}\geq \frac{n}{n-1}d_1\\
 e&\geq& \frac{d_2+\cdots +d_n}{n-2}\geq \frac{n-1}{n-2}d_2\\
\vdots\\
e&\geq& d_{n-1}+d_n \Rightarrow e\geq\frac{2}{1}d_{n-1}\\
\text{so,   } e^{n-1}&\geq&
\frac{n}{n-1}\cdot\frac{n-1}{n-2}\cdots\frac{2}{1}d_1d_2\cdots
d_{n-1}=n\prod_{i=1}^{n-1}d_i
\end{eqnarray*}
\begin{eqnarray*}
n!e(R/J)&=&n!(e^n-\prod_{i=1}^nd_i)=n!e^n-n!\prod_{i=1}^nd_i\geq
n!e^n-(n-1)!e^{n-1}d_n\\
&=&(n-1)e^{n-1}(ne-d_n)=\prod_{i=1}^nm_i.
\end{eqnarray*}
\end{proof}
\end{theorem}

\begin{lemma}\label{lemm2}
 If $M_1=e$, then $M_i=ie$ for
$i=1,...,n-1$. In particular $(k-1)e\leq\sum_{j=1}^kd_j$ for $2\leq k\leq n$.
\begin{proof}
Since $M_1=e$, $m_1=ne-(d_1+\cdots+d_n)$, $ne-(d_1+...+d_{n-1}) =ne-(d_1+...+d_n)+d_n\leq
e+d_n\leq 2e$. So, $M_2=2e$. $2e \geq ne-(d_1+ \cdots +d_{n-1})$, $ne-(d_1+...+d_{n-2}) =ne-(d_1+...+d_{n-1})+d_{n-1}\leq 2e+d_{n-1}\leq 3e$. So, $M_3=3e$. Repeat this process, $M_i=ie$ for all $i=1,...,n-1$ which implies automatically $m_i=ne-\sum_{j=i}^n d_j$ for $i=1,...,n$. $ie\geq ne-\sum_{j=1}^{n-i+1}d_j\geq ne-\sum_{j=i}^nd_j$. We obtain $(n-i)e\leq \sum_{j=1}^{n-i+1}d_j$ for $i=1,...,n-1$. So $(k-1)e\leq\sum_{j=1}^kd_j$ for $2\leq k\leq n$.
\end{proof}
\end{lemma}

\begin{theorem}
If  $M_1=e$, then the upper bound holds.
\begin{proof}
By lemma\  \ref{lemm2}, $M_i=ie$ for all $i=1,...,n-1$.
\begin{eqnarray*}
&e\geq ne-\sum_{j=1}^nd_j \Rightarrow e\leq \frac{1}{n-1}\sum_{j=1}^nd_j\leq\frac{n}{n-1}d_n,\\
&2e\geq ne-\sum_{j=1}^{n-1}d_j\Rightarrow e\leq \frac{1}{n-2}\sum_{j=1}^{n-1}d_j\leq\frac{n-1}{n-2}d_{n-1},\\
&\vdots\\
&(n-1)e\geq ne-(d_1+d_2)\Rightarrow e\leq d_1+d_2\leq 2d_2.
\end{eqnarray*}

We obtain
\begin{eqnarray*}
e^{n-1}\leq
\frac{n}{n-1}\cdot\frac{n-1}{n-2}\cdots\frac{2}{1}d_nd_{n-1}\cdots
d_2=n\prod_{i=2}^nd_i
\end{eqnarray*}

Now we show that $e(R/J)=e^n-\prod_{i=1}^nd_i\leq \frac{1}{n!}\prod_{i=1}^nM_i$.

\begin{eqnarray*}
\prod_{i=1}^nM_i&=&e(2e)\cdots
(n-1)e(ne-d_1)=n!e^n-(n-1)!e^{n-1}d_1\\
&\geq& n!e^n-n!\prod_{i=1}^nd_i,  \text{ because } e^{n-1}\leq
n\prod_{i=2}^nd_i
\end{eqnarray*}
\end{proof}
\end{theorem}


\end{document}